\documentclass[11pt]{article}
\usepackage{amsmath,amsfonts,amssymb,amsthm}
\title{Second moments related to \\Poisson hyperplane tessellations}
\author{Rolf Schneider\\
}
\date{}
\newtheorem{theorem}{Theorem}[section]

\newcommand{\Sd}{{\mathbb S}^{d-1}}
\newcommand{\Rd}{{\mathbb R}^d}
\newcommand{\D}{{\rm d}}

\newcommand{\E}{{\mathbb E}\,}
\newcommand{\eps}{\varepsilon}
\newcommand{\Hc}{\mathcal H}

\begin{document}

\maketitle

\begin{abstract}
It is well known that the vertex number of the typical $k$-face of a stationary random hyperplane tessellation in $\Rd$ has, under some mild conditions, an expectation equal to $2^k$, independent of the underlying distribution. The variance of this vertex number can vary widely. Under Poisson assumptions, we give sharp bounds for this variance, showing, in particular, that its maximum is attained if the tessellation is isotropic (that is, its distribution is rotation invariant) with respect to a suitable scalar product on $\Rd$.

The employed representation of the second moment of the vertex number is a special case of formulas providing the covariance matrix of the random vector $(\ell_0,\dots,\ell_k)$, where $\ell_r$ is the total $r$-face content of the typical $k$-face of a stationary Poisson hyperplane mosaic. For $k=d$ in the isotropic case, such formulas were first obtained by Miles in 1961. We give a more elementary proof and extend the formulas to general orientation distributions and to $k$-dimensional faces.
\\[2mm] 
\end{abstract}

\section{Introduction}\label{sec1}

Under suitable assumptions, a stationary random hyperplane process $\widehat X$ in $\Rd$ induces a tessellation of $\Rd$ into bounded polytopes. For $k\in\{1,\dots,d\}$, the typical $k$-face of such a tessellation, intuitively and heuristically speaking, is obtained by choosing at random, with equal chances, one of the $k$-dimensional polytopes of the tessellation within a `large' region of $\Rd$; a precise definition is recalled below. It is well known that a number of combinatorial quantities connected with this typical $k$-face have expectations that are essentially independent of the distribution of the underlying hyperplane process (see, e.g., Theorem 10.3.1 in \cite{SW08}). For example, the expected vertex number of the typical $k$-face is equal to $2^k$, which is the obvious value for a {\em parallel process}, where the hyperplanes belong to only $d$ translation classes. In the latter case, the typical $k$-face is a $k$-dimensional parallelepiped and thus its vertex number is constant, hence it has variance zero. In this note, we show, as a consequence of more general results, that the variance of the vertex number of the typical $k$-face of a stationary Poisson hyperplane mosaic attains its maximum in the isotropic case.

Let $\widehat X$ be a stationary Poisson hyperplane process in $\Rd$, with intensity measure $\widehat \Theta$ and intensity $\widehat \gamma>0$. Its spherical directional distribution is denoted by $\widehat\varphi$. We refer to \cite{SW08}, in particular Sections 4.4 and 10.3, for notions that are not explained here. We assume that $\widehat X$ is nondegenerate, which means that $\widehat \varphi$ is not concentrated on some great subsphere of $\Sd$. It follows from this assumption that $\widehat X$ almost surely induces a tessellation of $\Rd$ into bounded polytopes.  For $k=0, \dots,d$, we denote by ${\mathcal F}_k(\widehat X)$ the set of $k$-faces of the tessellation induced by $\widehat X$. As a particle process, this is stationary and hence has a grain distribution ${\mathbb Q}^{(k)}$, with respect to some center function. The {\em typical $k$-face} $Z^{(k)}$  of $X$ is defined as a random polytope with distribution ${\mathbb Q}^{(k)}$. (The choice of the center function is irrelevant, since we consider only translation invariant functions of the typical $k$-face.) In particular, the random system ${\mathcal F}_d(\widehat X)$ of $d$-dimensional polytopes or {\em cells} is the Poisson hyperplane mosaic induced by $\widehat X$. It is denoted by $X$ and its intensity by $\gamma$.

If $\Rd$ is equipped with a scalar product, then $\widehat X$ (as well as $X$) is called {\em isotropic} if its distribution is invariant under rotations. 

The vertex number $f_0=L_0$ is just the first polytope functional in the series $L_0,\dots,L_d$ of total $r$-face contents. These are defined by
\begin{equation}\label{1.1} 
L_r(P) =\sum_{F \in{\mathcal F}_r(P)} \Hc^r(F),
\end{equation}
for convex polytopes $P\subset \Rd$ and for $r\in\{0,\dots,d\}$. Here ${\mathcal F}_r(P)$ is the set of $r$-dimensional faces of $P$ and $\Hc^r$ denotes the $r$-dimensional Hausdorff measure. In particular, $L_d$ is the volume, $L_{d-1}$ the surface area, $L_1$ the total edge length, and $L_0=f_0$ the number of vertices. For the typical cell $Z$ of a stationary, isotropic Poisson hyperplane mosaic, Miles \cite{Mil61} has determined all mixed moments $\E(L_rL_s)(Z)$, $r,s\in\{0,\dots,d\}$. The result is reproduced, without proof, in \cite[formula (63)]{Mil70}. The proof given by Miles in \cite{Mil61} makes heavy use of ergodic theory and is not explicitly carried out in all details. Below, we give a short proof, based on the Slivnyak--Mecke formula, and extend the result to the typical $k$-face and to not necessarily isotropic mosaics (Theorem \ref{Thm2}). The latter is in contrast to a remark of Miles, who after the treatment of the isotropic case in \cite[Sec. 11.7]{Mil61} wrote: ``It does not seem at all practicable to generalise the theory of $\S\S$3--7 to the case of a general orientation distribution''. Of course, the general result cannot be so explicit as in the isotropic case, since the second moments heavily depend on the directional distribution of the underlying hyperplane process. The result makes this dependence as explicit as possible, in terms of the associated zonoid. 

The {\em associated zonoid} of the hyperplane process $\widehat X$ (see \cite[p. 156]{SW08}) is the convex body $\Pi_{\widehat X}$ with support function given by
\begin{equation}\label{2.10} 
h(\Pi_{\widehat X},u)= \frac{\widehat \gamma}{2} \int_{\Sd} |\langle u,v\rangle|\,\widehat\varphi(\D v),\quad u\in\Rd.
\end{equation}
Let $L\subset \Rd$ be a linear subspace of dimension at least $1$. Then the section process $\widehat X\cap L$ is a stationary Poisson hyperplane process with respect to $L$, and its associated zonoid is given by
$$ \Pi_{\widehat X\cap L} = \Pi_{\widehat X}|L,$$
where $(\cdot)|L$ denotes the orthogonal projection to $L$.

For a convex body $K$, we denote by $V_j(K)$ its $j$th intrinsic volume (see \cite[Sec. 14.2]{SW08} or \cite[Chap. 4]{Sch14}). If $K$ is centrally symmetric with respect to $0$ (as it holds for $\Pi_{\widehat X}$ and its projections), then $K^\circ$ denotes the polar body of $K$, constructed within the linear hull of $K$.

In the subsequent theorems, $\nabla_j(u_1,\dots,u_j)$ denotes the $j$-dimensional volume of the parallelepiped spanned by the vectors $u_1,\dots,u_j$, and $\kappa_j$ is the $j$-dimensional volume of the $j$-dimensional unit ball.

\begin{theorem}\label{Thm2}
The total face contents of the typical $k$-face $Z^{(k)}$ of the stationary Poisson random mosaic $X$ satisfy
\begin{eqnarray}\label{4.10}
&& \E (L_rL_s)(Z^{(k)}) \nonumber\\
&& = \sum_{j=\max\{r,s\}}^{k} \frac{k!j!}{(k-j)!(j-s)!}\, 2^{k-2j} \nonumber\\
&& \hspace{4mm}\times\; \frac{\widehat\gamma^{d-s}}{\gamma d!} \int_{(\Sd)^{d-s}} V_{j-r}(\Pi_{\widehat X}|u_1^\perp\cap\dots\cap u_{d-j}^\perp)V_j((\Pi_{\widehat X}|u_1^\perp\cap\dots\cap u_{d-j}^\perp)^\circ)\nonumber\\
&&  \hspace{4mm}\times\;\nabla_{d-s}(u_1,\dots,u_{d-s})\,\widehat\varphi^{\hspace{1pt}{d-s}}(\D(u_1,\dots,u_{d-s}))
\end{eqnarray} 
for $k\in\{1,\dots,d\}$ and $r,s\in\{0,\dots,d\}$.
\end{theorem}

Together with the known relations
\begin{equation}\label{4.18} 
\E L_r(Z^{(k)}) =\frac{2^{k-r}\binom{k}{r}}{\gamma\binom{d}{r}}\,V_{d-r}(\Pi_{\widehat X}), \quad r=0,\dots,d,
\end{equation}
following from (\ref{3.5}) and (\ref{3.6}) below, Theorem \ref{Thm2} provides an explicit expression for the covariance matrix of the random vector $(L_0(Z^{(k)}),\dots,L_k(Z^{(k)}))$. 

For $r=0$, the functional in the integrand of (\ref{4.10}) contains a so-called volume product, for which (in the cases appearing here) sharp bounds are known. The remaining integral, together with its factor, can then be evaluated. In this way, we obtain the following result, which can be considered as a reformulation of a result obtained in \cite{Sch10}.

\begin{theorem}\label{Thm1} 
Let $\widehat X$ be a nondegenerate stationary Poisson hyperplane process in $\Rd$ $($$d\ge 2$$)$, and let $Z^{(k)}$ be the typical $k$-face of its induced mosaic $($$k\in\{2,\dots,d\}$$)$. The variance of the vertex number $f_0(Z^{(k)})$ satisfies
\begin{equation}\label{1.0} 
0 \le {\rm Var}f_0(Z^{(k)}) \le 2^k k!\left(\sum_{j=0}^k \frac{\kappa_j^2}{4^j(k-j)!}\right)  -2^{2k}.
\end{equation}
Equality on the left side of $(\ref{1.0})$ holds if and only if $\widehat X$ is a parallel process. Equality on the right side holds if  $\widehat X$ is isotropic with respect to a suitable scalar product on $\Rd$, and for $k=d$ it holds only in this case.
\end{theorem}

We remark that also for the volume $V_d=L_d$, estimates of ${\rm Var}\,V_d(Z)$ can be obtained. But since $V_d(Z)$ is not homogeneous of degree $0$, no absolute inequalities are possible, but only inequalities of isoperimetric type. It follows from (\ref{4.0}), (\ref{3.7}), (\ref{4.18}) below that for each stationary Poisson hyperplane process $\widehat X$ we have
$$ \frac{{\rm Var}\,V_d(Z)}{({\mathbb E}V_d(Z))^2} = 2^{-d}d!{\rm vp}(\Pi_{\widehat X})-1.$$
Hence, by (\ref{5.2}) and its equality cases, for given expectation ${\mathbb E}V_d(Z)$, the variance ${\rm Var}\,V_d(Z)$ is maximal for isotropic processes and minimal for parallel processes.

In the next section, after some preparations, we recall a result of Favis and Weiss \cite{FW98} on $r$-face-content weighted typical cells and its extension to weighted typical $k$-faces. This is used in Section \ref{sec3} to prove Theorem \ref{Thm2}. The proof of Theorem \ref{Thm1} is then explained in Section \ref{sec4}. That section also mentions some stability results in the case $k=d$.

\section{Weighted faces}\label{sec2}

We fix some notation and recall some definitions. We work in the $d$-dimensional real vector space $\Rd$ ($d\ge 2$) and use its standard scalar product $\langle\cdot,\cdot\rangle$ to define, for example, its unit sphere $\Sd$. The space of hyperplanes in $\Rd$ is equipped with its usual topology. The hyperplane through $0$ orthogonal to $u\in \Sd$ is denoted by $u^\perp$.

By ${\mathcal P}^d$ we denote the space of (nonempty, compact, convex) polytopes in $\Rd$, endowed with the topology induced by the Hausdorff metric. For a polytope $P\in{\mathcal P}^d$ and for $r\in\{0,\dots,d\}$, we denote, as already mentioned, by ${\mathcal F}_r(P)$ the set of its $r$-dimensional faces.

If ${\mathcal H}$ is any locally finite system of hyperplanes in ${\mathbb R}^d$ that induces a tessellation of ${\mathbb R}^d$ into bounded polytopes, then we denote by ${\mathcal F}_k({\mathcal H})$ the set of all $k$-dimensional polytopes of this tessellation (that is, the set of all $k$-faces of its $d$-polytopes).

Probabilities are denoted by ${\mathbb P}$ and expectations by $\E$. For a topological space $T$, we denote by ${\mathcal B}(T)$ the $\sigma$-algebra of its Borel sets. When we use simple counting measures below (implicitly in the use of point processes), it is convenient to identify a simple counting measure $\eta$ with its support; with this identification made, $\eta(\{x\})=1$ and $x\in\eta$ are used synonymously, and so are $\eta \cup \eta'$ and $\eta + \eta'$, as long as the latter is simple. 

We recall from \cite{Sch10} some basic facts about weighted faces. Let $k\in\{1,\dots,d\}$. Let $w$ be a positive, translation-invariant, measurable function on the space ${\mathcal P}^{(k)}$ of $k$-dimensional convex polytopes in ${\mathbb R}^d$, with finite expectation ${\mathbb E} \,w(Z^{(k)})$. The {\em $w$-weighted typical $k$-face}  $Z_w^{(k)}$ of $X$ is defined as a random polytope with distribution given by 
$$ {\mathbb P}(Z^{(k)}_w\in A) = \frac{1}{{\mathbb E} w(Z^{(k)})}\, {\mathbb E} \left[{\bf 1}_A(Z^{(k)})w(Z^{(k)})\right]$$
for $A\in{\cal B}({\mathcal P}^{(k)})$; here ${\bf 1}_A$ denotes the indicator function of $A$. For a nonnegative, translation-invariant, measurable function $h$ on ${\mathcal P}^{(k)}$ we have
$${\mathbb E}h(Z^{(k)}_w) = \lim_{r\to\infty} \frac{{\mathbb E} \sum_{K\in {\mathcal F}_k(\widehat X),\,K\subset rW} h(K) w(K)}{{\mathbb E}\sum_{K\in {\mathcal F}_k(\widehat X),\,K\subset rW} w(K)} $$
for $A\in{\cal B}({\mathcal P}^{(k)})$ and any compact convex set $W$ of positive volume, as follows from \cite[Thm. 4.1.3]{SW08}. This gives some idea of the intuitive meaning of the weighted typical $k$-face, and it also shows that
\begin{equation}\label{3.1}
{\mathbb E}(hw)(Z^{(k)}) = {\mathbb E}h(Z_w^{(k)})\cdot{\mathbb E}w(Z^{(k)}).
\end{equation}

As weights we now use the functions $L_s$ defined by (\ref{1.1}). A formula for the distribution (up to translations) of the $L_s$-weighted typical cell was derived by Favis and Weiss \cite{FW98}. With a slightly simplified proof, using the Slivnyak--Mecke formula, this was extended in \cite[Thm. 1]{Sch10} to $L_s$-weighted typical $k$-faces. For $k=1,\dots,d$, $s=0,\dots,k$ and with $h$ as above, the result can be formulated as 
\begin{eqnarray}\label{3.2}
\E h(Z^{(k)}_{L_s}) &=& \frac{\widehat\gamma^{d-s}}{\widehat\gamma_{d-s}}\frac{2^{s-k}}{(d-s)! \binom{d-s}{d-k}}
\int_{(\Sd)^{d-s}} \E\sum_{0\in K\in{\mathcal F}_k(\widehat X\cup\{u_1^\perp,\dots, u_{d-s}^\perp\})} h(K)\nonumber\\
&& \times\,\nabla_{d-s}(u_1,\dots,u_{d-s})\,\widehat\varphi^{\hspace{1pt}d-s}(\D(u_1,\dots,u_{d-s})).
\end{eqnarray}
We have rewritten the result, using definition (5) of \cite{Sch10} (and corrected a misprint). Here, $\widehat\gamma_{d-s}$ is the intensity of the intersection process of order $d-s$ of the hyperplane process $\widehat X$ (see \cite[pp. 133--134]{SW08}).

To further rewrite (\ref{3.2}), we state that
\begin{equation}\label{3.3}
\frac{{\mathbb E}L_s(Z^{(k)})}{\widehat\gamma_{d-s}} = \frac{2^{k-s}\binom{k}{s}}{\binom{d}{s}\gamma}.
\end{equation}
Indeed, by \cite[(4.63)]{SW08} we have
\begin{equation}\label{3.4}
\widehat\gamma_{d-s} =V_{d-s}(\Pi_{\widehat X}).
\end{equation}
From \cite[(7)]{Sch13}, for $r=s$, it follows that
\begin{equation}\label{3.5}
{\mathbb E}L_s(Z^{(k)}) = 2^{k-s}\binom{k}{s}{\mathbb E}V_s(Z^{(s)}),
\end{equation}
and from \cite[(12)]{Sch13} we get
\begin{equation}\label{3.6}
{\mathbb E}V_s(Z^{(s)}) =\frac{V_{d-s}(\Pi_{\widehat X})}{\binom{d}{s}\gamma},
\end{equation}
where we have used that
\begin{equation}\label{3.7} 
\gamma = V_d(\Pi_{\widehat X}).
\end{equation}
This is formula (10.44) in \cite{SW08} for $k=d$ (note that $\gamma^{(d)}$ appearing there is the intensity of the process of $d$-dimensional cells of $X$ and hence is what we have here denoted by $\gamma$). Relations (\ref{3.4})--(\ref{3.6}) together yield (\ref{3.3}).

From (\ref{3.1})--(\ref{3.3}) we now obtain that
\begin{eqnarray}\label{3.8}
\E (hL_s)(Z^{(k)}) &=& \frac{\widehat\gamma^{d-s}}{\gamma} \frac{1}{(d-s)!\binom{d}{k}} 
\int_{(\Sd)^{d-s}} \E\sum_{0\in K\in{\mathcal F}_k(\widehat X\cup\{u_1^\perp,\dots, u_{d-s}^\perp\})} h(K)\nonumber\\
&& \times\,\nabla_{d-s}(u_1,\dots,u_{d-s})\,\widehat\varphi^{\hspace{1pt}d-s}(\D(u_1,\dots,u_{d-s})).
\end{eqnarray}

If $s=d$ (and hence $k=d$), no additional hyperplanes $u_1^\perp,\dots,u_{d-s}^\perp$ occur, and this formula has to be read with the integrations deleted, thus
\begin{equation}\label{3.9} 
{\mathbb E}(hL_d)(Z^{(d)})=\frac{1}{\gamma} \,{\mathbb E}\sum_{0\in K\in{\mathcal F}_d(\widehat X)} h(K) =\frac{1}{\gamma} \,{\mathbb E}h(Z_0),
\end{equation}
where $Z_0$ denotes the zero cell of $X$, that is, the almost surely unique cell of $X$ that contains the origin. This just expresses the well-known fact that the distribution of the  zero cell, if translations are disregarded, is the volume-weighted distribution of the typical cell (\cite[Thm. 10.4.1]{SW08}).

\section{A covariance matrix}\label{sec3}

In this section, we determine $\E(L_rL_s)(Z^{(k)})$ for $k\in\{1,\dots,d\}$ and thus prove Theorem \ref{Thm2}. First we consider the case $s=d$ and hence $k=d$. Relation (\ref{3.9}) with $h=L_r$, together with \cite[(10.51)]{SW08}, yields
\begin{equation}\label{4.0}
\E(L_rL_d)(Z) = \frac{d!}{2^d \gamma} V_{d-r}(\Pi_{\widehat X})V_d(\Pi_{\widehat X}^\circ).
\end{equation}

Now let $s \in\{0,\dots,d-1\}$ and $r\in\{0,\dots,d\}$ with $r,s\le k$. We want to apply (\ref{3.8}) with $h=L_r$. To transform the sum $\sum h(K)$ in (\ref{3.8}), which extends over the $k$-faces of the tessellation induced by $\widehat X\cup\{u_1^\perp, \dots, u_{d-s}^\perp\}$ and containing $0$, we generalize an idea of Miles \cite[Sec. 11.6]{Mil61}.

Let $H_1,\dots,H_{d-s}$ be fixed hyperplanes through $0$ in general position (hyperplanes in ${\mathbb R}^d$ are said to be {\em in general position} if any $m\le d$ of them have an intersection of dimension $d-m$). Almost surely, $\widehat X$ and $H_1,\dots,H_{d-s}$ are in general position, and no hyperplane of $\widehat X$ passes through $0$; this is assumed for the realisations of $\widehat X$ considered in the following. We define
$$ {\mathcal C}_k := \big\{ K\in {\mathcal F}_k(\widehat X\cup\{H_1,\dots,H_{d-s}\}): 0\in K\big\}.$$

To evaluate (\ref{3.8}) for $h=L_r$, we have to determine
$$ \sum_{K\in{\mathcal C}_k} L_r(K) = \sum_{K\in{\mathcal C}_k}\; \sum_{F\in{\mathcal F}_r(K)} \Hc^r(F).$$
Let $F\in{\mathcal F}_r(K)$ for some $K\in{\mathcal C}_k$. Since $K$ contains $0$, it is not contained in a hyperplane of $\widehat X$. Therefore, there are $d-k$ hyperplanes $H_{l_1},\dots,H_{l_{d-k}}\in\{H_1,\dots,H_{d-s}\}$ such that
$$ K\subset Z_0\cap H_{l_1}\cap\dots\cap H_{l_{d-k}}.$$
The $r$-face $F$ is the intersection of $K$ with $k-r$ hyperplanes from $\widehat X\cup \{H_1,\dots,H_{d-s}\}\setminus \{H_{i_1},\dots, H_{i_{d-k}}\}$. Suppose that the hyperplanes from $\{H_1,\dots,H_{d-s}\}$ that contain $F$ have an intersection of dimension $j$. Then $\max\{r,s\}\le j\le k$.

Define 
$$ {\mathcal Z}_j := \big\{ Z_0\cap H_{m_1}\cap\dots\cap H_{m_{d-j}}: 1\le m_1<\dots< m_{d-j}\le d-s\big\}$$
for $j=\max\{r,s\},\dots,k$ (with ${\mathcal Z}_d:=\{Z_0\}$ if $k=d$). Thus, ${\mathcal Z}_j$ is a set of $j$-dimensional polytopes, each containing the origin. The (generalized) approach of Miles consists in replacing the summation over the $r$-faces of the polytopes in ${\mathcal C}_k$ by a summation over the $r$-faces of the polytopes in the families ${\mathcal Z}_j$, regarding multiplicities. This has the advantage that $\E L_r(Z_0\cap L)$, for a linear subspace $L$, can be evaluated. First we determine the multiplicities.

Let $P:=Z_0 \cap H_{m_{1}}\cap\dots\cap H_{m_{d-j}}\in {\mathcal Z}_j$ (with $1\le m_1<\dots< m_{d-j}\le d-s$ and $\max\{r,s\} \le j\le k$) be given, and let $G$ be one of its $r$-faces. We assume first that $j>0$. By $H^+,H^-$ we denote the two closed halfspaces bounded by a hyperplane $H$ through $0$. Let $K\in{\mathcal C}_k$ be a polytope that has an $r$-face $F$ contained in $G$. From $F\subset G\subset P$ it follows that $F\subset H_{m_1},\dots,H_{m_{d-j}}$. Since $j\not=0$, the face $G$ is the intersection of $H_{m_{1}}\cap\dots\cap H_{m_{d-j}}$ with a $(d+r-j)$-face of $Z_0$. Since any intersection with a further hyperplane from $\{H_1,\dots,H_{d-s}\}\setminus \{H_{m_1},\dots, H_{m_{d-j}}\}$ reduces the dimension (by general position), it follows that $F\not\subset H_m$ for $m\in\{1,\dots,d-s\}\setminus\{m_1,\dots,m_{d-j}\}$. Therefore, we can write
\begin{equation}\label{4.1} 
K= Z_0\cap H_{l_1}\cap\dots\cap H_{l_{d-k}}\cap \left(H^{\varepsilon_{d-k+1}}_{l_{d-k+1}}\cap\dots \cap H^{\varepsilon_{d-j}}_{l_{d-j}}\right)
\cap H^{\varepsilon_{d-j+1}}_{l_{d-j+1}}\cap\dots \cap H^{\varepsilon_{d-s}}_{l_{d-s}},
\end{equation}
where 
$$ l_1,\dots,l_{d-k}\in \{m_1,\dots,m_{d-j}\},$$
further $\{l_1,\dots,l_{d-s}\}=\{1,\dots,d-s\}$, and $\varepsilon_m\in\{+,-\}$ for $m=d-k+1,\dots,d-s$. Here the indices are chosen such that each of the hyperplanes $H_{l_1}, \dots, H_{l_{d-k}}$ contains $K$, each of the hyperplanes $H_{l_{d-k+1}}, \dots, H_{l_{d-j}}$ contains $F$, and each of the hyperplanes $H_{l_{d-j+1}},\dots,H_{l_{d-s}}$ does not contain $F$. If in (\ref{4.1}) we change any of the signs $\varepsilon_{d-k+1}, \dots,\varepsilon_{d-j}$ (but fix $\varepsilon_{d-j+1}, \dots,\varepsilon_{d-s}$), then the resulting polytope $K$ still belongs to ${\mathcal C}_k$ and contains $F$. All polytopes $K\in{\mathcal C}_k$ containing $F$ are obtained in this way. It follows that $F$ belongs to precisely $\binom{d-j}{k-j}2^{k-j}$ polytopes from ${\mathcal C}_k$. In the case $j=0$, which was excluded up to now, a simpler variant of the preceding argument shows that $\{0\}$  belongs to precisely $\binom{d}{k}2^{k}$ polytopes from ${\mathcal C}_k$. Moreover, if $j,r>0$, each $r$-face $G$ of a polytope $P$ as above is divided by hyperplanes from $H_1,\dots,H_{d-s}$ into $r$-faces of polytopes from ${\mathcal C}_k$, so that the $r$-volumes of the latter add up to the $r$-volume of $G$. We conclude that
\begin{eqnarray*}
\sum_{K\in{\mathcal C}_k} L_r(K) &=& \sum_{K\in{\mathcal C}_k}\; \sum_{F\in{\mathcal F}_r(K)} \Hc^r(F)\\
&=& \sum_{j=\max\{r,s\}}^k \binom{d-j}{k-j}2^{k-j} \sum_{P\in{\mathcal Z}_j} \,\sum_{G\in{\mathcal F}_r(P)} \Hc^r(G)\\
&=& \sum_{j=\max\{r,s\}}^k \binom{d-j}{k-j}2^{k-j} \sum_{P\in{\mathcal Z}_j} L_r(P).
\end{eqnarray*} 

Inserting this into (\ref{3.8}) with $h=L_r$, we obtain
\begin{eqnarray}\label{4n3}
&& \E (L_rL_s)(Z^{(k)}) \nonumber\\
&& =\frac{\widehat\gamma^{d-s}}{\gamma (d-s)!} \binom{d}{k}^{-1} \sum_{j=\max\{r,s\}}^{k} \binom{d-j}{k-j} 2^{k-j}
\int_{(\Sd)^{d-s}} \E \sum_{1 \le m_1 < \dots < m_{d-j} \le d-s} \nonumber\\
&& \hspace{4mm}\times\; L_r(Z_0\cap u_{m_1}^\perp \cap \dots \cap u_{m_{d-j}}^\perp)\nabla_{d-s}(u_1,\dots,u_{d-s})\,\widehat\varphi^{\hspace{1pt}{d-s}}(\D(u_1,\dots,u_{d-s}))\nonumber\\
&& =\sum_{j=\max\{r,s\}}^{k} \frac{k!}{(k-j)!(j-s)!}\, 2^{k-j} \frac{\widehat\gamma^{d-s}}{\gamma d!} \int_{(\Sd)^{d-s}} \E L_r(Z_0\cap u_1^\perp \cap \dots \cap u_{d-j}^\perp)\nonumber\\
&&  \hspace{4mm}\times\;\nabla_{d-s}(u_1,\dots,u_{d-s})\,\widehat\varphi^{\hspace{1pt}{d-s}}(\D(u_1,\dots,u_{d-s})).
\end{eqnarray}

For given linearly independent vectors $u_1,\dots,u_{d-j}\in \Sd$ in general position, let
$$ L:= u_1^\perp\cap\dots\cap u_{d-j}^\perp.$$
We can now argue as in \cite[p. 690]{Sch09}: the intersection $Z_0\cap L$ is the zero cell of the intersection process $\widehat X\cap L$ (see \cite[(4.61)]{SW08}), and it is known (see \cite[Thm. 10.4.9]{SW08}) that
$$ \E L_r(Z_0\cap L) = 2^{-j}j! V_{j-r}(\Pi_{\widehat X}|L)V_j((\Pi_{\widehat X}|L)^\circ).$$
Thus, the expectation in the integrand of (\ref{4n3}) can be expressed in terms of intrinsic volumes of the associated zonoid.  This completes the proof of Theorem \ref{Thm2}.

Now we specialize Theorem \ref{Thm2} to the isotropic case. Assume that $\widehat X$ is isotropic. Then the associated zonoid $\Pi_{\widehat X}$ is a ball of radius
\begin{equation}\label{4.15} 
R= \widehat\gamma \, \frac{\kappa_{d-1}}{d\kappa_d}, 
\end{equation}
see \cite[p. 490]{SW08}. The $k$th intrinsic volume of the unit ball $B^d$ is given by
$$ V_k(B^d) =\binom{d}{k}\frac{\kappa_d}{\kappa_{d-k}}$$
(\cite[(14.8)]{SW08}). It follows that 
\begin{eqnarray*} 
V_{j-r}(\Pi_{\widehat X}|u_1^\perp\cap\dots\cap u_{d-j}^\perp) &=& R^{j-r} \binom{j}{r}\frac{\kappa_j}{\kappa_r},\\
V_j((\Pi_{\widehat X}|u_1^\perp\cap\dots\cap u_{d-j}^\perp)^\circ) &=& R^{-j}\kappa_j.
\end{eqnarray*}
The remaining integral is known (also in the non-isotropic case), namely
\begin{equation}\label{4.17}  
\widehat\gamma^{\hspace{1pt}d-s} \int_{(\Sd)^{d-s}} \nabla_{d-s}(u_1,\dots,u_{d-s})\,\widehat\varphi^{\hspace{1pt}d-s}(\D(u_1,\dots,u_{d-s}))=(d-s)!V_{d-s}(\Pi_{\widehat X}),
\end{equation}
by \cite[(14.35)]{SW08} with $\rho=\widehat\gamma\widehat\varphi/2$. Since 
$$ V_{d-s}(\Pi_{\widehat X}) = R^{d-s}\binom{d}{s}\frac{\kappa_d}{\kappa_s},\qquad V_d(\Pi_{\widehat X})=R^d\kappa_d,$$
we obtain from Theorem \ref{Thm2} that
\begin{equation}\label{4.16}  
\E(L_rL_s)(Z^{(k)}) = \frac{2^k k!}{\kappa_r\kappa_s} \left(\frac{d\kappa_d}{\kappa_{d-1}\widehat\gamma}\right)^{r+s}  \sum_{j=\max\{r,s\}}^k \frac{\kappa_j^2}{4^j(k-j)!}\binom{j}{r}\binom{j}{s}.
\end{equation}
Since we can also write
$$ \frac{d\kappa_d}{\kappa_{d-1}} = \frac{2\pi^{\frac{1}{2}} \Gamma \left(\frac{1}{2} \left[d+1\right] \right)}
{\Gamma\left(\frac{1}{2}d\right)},\qquad \kappa_j = \frac{\pi^{\frac{j}{2}}}{\Gamma\left(\frac{j}{2}+1\right)} = \frac{2^j\pi^{\frac{j-1}{2}}\Gamma\left(\frac{1}{2}\left[j+1\right]\right)}{j!}, $$
(\ref{4.16}) is the same as
\begin{eqnarray*}
\E(L_rL_s)(Z^{(k)}) &=& \frac{2^k \pi^{\frac{1}{2}}}{\Gamma\left(\frac{1}{2}\left[r+1\right]\right) \Gamma\left(\frac{1}{2}\left[s+1\right]\right)} \left\{\frac{\Gamma\left(\frac{1}{2}\left[d+1\right]\right)} {\Gamma\left(\frac{1}{2}d\right)\widehat\gamma}\right\}^{r+s}\\
&& \times \sum_{j=\max\{r,s\}}^k \binom{k}{j} \left(\frac{\pi}{2}\right)^j \frac{\Gamma \left(\frac{1}{2} \left[j+1\right]\right)}{\Gamma\left(\frac{1}{2}j+1\right)}(j)_r(j)_s
\end{eqnarray*}
with $(j)_r=j!/(j-r)!$. For $k=d$, this is formula (63) of Miles \cite{Mil70}.

Another simple case where (\ref{4.10}) can be evaluated explicitly is that of a cuboid process $\widehat X$, where the hyperplanes of the process belong to $d$ pairwise orthogonal translation classes. Cuboid processes were studied by Favis \cite{Fav96}. We consider only the case $k=d$ for a quasi-isotropic cuboid process, where the directional distribution $\widehat\varphi$ is concentrated in $\pm e_1,\dots,\pm e_n$ for an orthonormal basis $(e_1,\dots,e_n)$ of $\Rd$ and assigns the same value to each of these points. In this case, the associated zonoid $\widehat X$ is a cube of edge length $2^{-d}\,\widehat\gamma$. Moreover, for linearly independent vectors $u_1,\dots,u_{d-j}$ in the support of the measure $\widehat\varphi$, the orthogonal projection $\Pi_{\widehat X}|u_1^\perp\cap\dots\cap u_{d-j}^\perp$ is a $j$-dimensional cube of edge length $2^{-d}\,\widehat\gamma$. For a $d$-dimensional cube $C$ of edge length $a$ and for $j\in\{0,\dots,d\}$ we have
$$ V_j(C)=\binom{d}{j}a^j$$
by \cite[(4.23)]{Sch14}, together with $f_j(C)=2^{d-j}\binom{d}{j}$. From the well-known relation $V_d(C)V_d(C^\circ)=4^d/d!$, if $C$ is centred at $0$, it follows that 
$$ V_d(C^\circ)= \frac{4^d}{d!a^d}.$$
Now we conclude from (\ref{4.10}), observing (\ref{4.17}), that
\begin{equation}
\E(L_rL_s)(Z)= \frac{2^{d(r+s+1)}}{\widehat\gamma^{\hspace{1pt}r+s}}  \sum_{j=\max\{r,s\}}^d \binom{d}{j} \binom{j}{r}\binom{j}{s}
\end{equation}
for $r,s\in\{0,\dots,d\}$.

\section{Proof of Theorem \ref{Thm1}}\label{sec4}

Essentially, Theorem \ref{Thm1} was already proved in \cite{Sch10}. More precisely, it follows from the inequalities, established there, for the expected vertex number of the $L_j$-weighted typical $k$-face. For the reader's convenience, we recall the final argument in our present case.

For $r=0$, the functional appearing in the integrand of (\ref{4.10}) is a volume product. For a $0$-symmetric convex body $K\subset\Rd$ of dimension $j\in \{0,\dots,d\}$, the {\em volume product} is defined by
$$ {\rm vp}(K):= V_j(K)V_j(K^\circ),$$
where the polar $K^\circ$ is taken with respect to the linear subspace spanned by $K$. Specializing (\ref{4.10}) to $r=s=0$, we obtain
\begin{eqnarray}\label{5.1}
\E f_0^2(Z^{(k)}) &=& \sum_{j=0}^k \frac{k!}{(k-j)!}2^{k-2j} \\
&& \times\, \frac{\widehat\gamma^d}{\gamma\, d!}\int_{(\Sd)^d}
{\rm vp}(\Pi_{\widehat X}|u_1^\perp\cap\dots\cap u_{d-j}^\perp)\nabla_d(u_1,\dots,u_d) \,\widehat \varphi^{\hspace{1pt}d}(\D(u_1,\dots,u_d)).\nonumber
\end{eqnarray}
We notice that for the integral appearing here we have
$$ \frac{\widehat\gamma^d}{\gamma\, d!}\int_{(\Sd)^d} \nabla_d(u_1,\dots,u_d) \,\widehat \varphi^{\hspace{1pt}d} (\D(u_1,\dots,u_d))=1.$$
This follows from (\ref{4.17}) for $s=0$, together with (\ref{3.7}).

For a $j$-dimensional zonoid $K$, such as $\Pi_{\widehat X}| u_1^\perp\cap\dots\cap u_{d-j}^\perp$, the inequalities
\begin{equation}\label{5.1a} 
\frac{4^j}{j!} \le {\rm vp}(K) \le \kappa_j^2
\end{equation}
are valid. The right-hand inequality is known as the Blaschke--Santal\'{o} inequality and the left-hand side (valid for zonoids) as Reisner's inequality (for references, see \cite[Chap. 14]{SW08}). We conclude that
$$
2^{2k} \le \E f_0^2(Z^{(k)}) \le \sum_{j=0}^k \frac{k!}{(k-j)!}2^{k-2j} \kappa_j^2.
$$
Since $\E f_0(Z^{(k)})=2^k$, this yields the inequalities (\ref{1.0}).´

We discuss the equality cases. If $\widehat X$ is a parallel process, then the typical $k$-face is a parallelepiped, hence $f_0(Z^{(k)})=2^k$ and thus ${\rm Var} f_0(Z^{(k)})=0$. Conversely, if equality holds on the left side of (\ref{1.0}), then it holds on the left side of (\ref{5.1a}), and hence $K$ is a parallelepiped, if $K= \Pi_{\widehat X}|u_1^\perp\cap\dots\cap u_{d-j}^\perp$ and $u_1,\dots,u_{d-j}$ are linearly independent vectors in the support of the measure $\widehat\varphi$. It was proved in \cite{Sch09} that then $\Pi_{\widehat X}$ is a parallelepiped and hence $\widehat X$ is a parallel process. If $\widehat X$ is isotropic with respect to some scalar product, then $\Pi_{\widehat X}$ is an ellipsoid, hence each $K=\Pi_{\widehat X}|u_1^\perp\cap\dots\cap u_{d-j}^\perp$ with linearly independent vectors $u_1,\dots,u_{d-j}$ is a $j$-dimensional ellipsoid, and equality holds on the right side of (\ref{5.1a}) and hence in (\ref{1.0}). The converse we can only show for $k=d$. In that case, equality on the right side of (\ref{1.0}) implies equality on the right inequality of
\begin{equation}\label{5.2}
\frac{4^d}{d!} \le {\rm vp}(\Pi_{\widehat X})\le \kappa_d^2.
\end{equation}
Therefore, $\Pi_{\widehat X}$ is an ellipsoid. Choosing a suitable scalar product on $\Rd$, we can assume that $\Pi_{\widehat X}$ is a ball. Since the associated zonoid $\Pi_{\widehat X}$ determines the distribution of $\widehat X$ uniquely (\cite[Thm. 4.6.4]{SW08}), $\widehat X$ is isotropic with respect to the new scalar product. 

Finally, we mention that the just proved characterizations for $k=d$, that is, the characterizations of the processes for which ${\rm Var}\,f_0(Z)$, the variance of the typical cell, attains its extreme values, can be strengthened in the form of stability assertions. Here we follow the example set out by B\"{o}r\"{o}czky and Hug \cite{BH10}. For this, we write
$$ {\rm Var}\,f_0(Z) = 2^{-d}d!\left[{\rm vp}(\Pi_{\widehat X})+\Phi(\widehat X)\right]$$
with
\begin{eqnarray*}
\Phi(\widehat X) &= & \sum_{j=0}^{d-1} \frac{2^{2(d-j)}}{(d-j)!} \cdot \frac{\widehat\gamma^d}{\gamma\, d!}\int_{(\Sd)^d}\\
&& \times\;{\rm vp}(\Pi_{\widehat X}|u_1^\perp\cap\dots\cap u_{d-j}^\perp)\nabla_d(u_1,\dots,u_d) \,\widehat\varphi^{\hspace{1pt}d} (\D(u_1,\dots,u_d)) -\frac{2^{3d}}{d!}.
\end{eqnarray*}
Then (\ref{5.2}) holds together with
$$ -\frac{4^d}{d!}=:c_d\le \Phi(\widehat X)\le C^d := \sum_{j=0}^{d-1} \frac{2^{2(d-j)}}{(d-j)!}\,\kappa_j^2 - \frac{2^{3d}}{d!}.$$

Suppose that ${\rm Var}\,f_0(Z)$ is close to its minimal value $0$, say
$$ {\rm Var}\,f_0(Z) \le \eps$$
with some $\eps>0$. Then 
$$ 2^{-d}d![{\rm vp}(\Pi_{\widehat X})+c_d] \le  2^{-d}d![{\rm vp}(\Pi_{\widehat X})+\Phi(\widehat X)] ={\rm Var}\,f_0(Z)\le \eps$$
and hence
\begin{equation}\label{5.3} 
{\rm vp}(\Pi_{\widehat X}) \le (1+2^{-d}\eps)\frac{4^d}{d!}.
\end{equation}
If (\ref{5.3}) holds with sufficiently small $\eps>0$, then it is shown in \cite{BH10} how a suitable distance (e.g., Wasserstein or Prokhorov) of the directional distribution $\widehat\varphi$ of $\widehat X$ from the distribution of a suitable parallel process can be estimated from above in terms of $\eps$.

Now suppose that ${\rm Var}\,f_0(Z)$ is close to its maximal value, say
$$ {\rm Var}\,f_0(Z) \ge (1-\eps)2^{-d}d![\kappa_d^2+C_d]$$
with some $\eps\in(0,1)$. Then 
$$ 2^{-d}d![{\rm vp}(\Pi_{\widehat X})+C_d] \ge  2^{-d}d![{\rm vp}(\Pi_{\widehat X})+\Phi(\widehat X)] ={\rm Var}\,f_0(Z)\ge (1-\eps) 2^{-d}d!(\kappa_d^2+C_d)$$
and hence
\begin{equation}\label{5.4} 
{\rm vp}(\Pi_{\widehat X}) \ge (1-a_d\eps)\kappa_d^2
\end{equation}
with 
$$ a_d = 1+C_d/\kappa_d^2>0.$$ 
If (\ref{5.4}) holds with sufficiently small $\eps>0$, then one can see from \cite{BH10} how a suitable distance of the directional distribution $\widehat\varphi$ of $\widehat X$ from the isotropic distribution can be estimated from above in terms of $\eps$.

\vspace{5mm}

\noindent Author's address:\\[2mm]
Rolf Schneider\\
Mathematisches Institut, Albert-Ludwigs-Universit{\"a}t\\
D-79104 Freiburg i. Br., Germany\\
E-mail: rolf.schneider@math.uni-freiburg.de

\end{document}